\documentclass[12pt]{article}

{\vspace{1em}
\begin{minipage}{12cm}
\rule{12cm}{1pt}\sffamily\noindent}%
{\newline\rule{12cm}{1pt}
\end{minipage}
 \vspace{1em}}

\usepackage{amsmath}
\usepackage{amsfonts}
\usepackage{amssymb}
\usepackage{latexsym}
\usepackage{euscript}
\usepackage{xspace}
\usepackage{gastex}

\newtheorem{teo}{Theorem}
\newtheorem{lem}[teo]{Lemma}
\newtheorem{pro}[teo]{Proposition}
\newtheorem{cor}[teo]{Corollary}

\newenvironment{proof}{%
{\noindent\em Proof}:}{\hspace*{\fill}$\Box$\vspace{3 mm}}
\newtheorem{problem}{Problem}

{\begin{enumerate}}%
{\end{enumerate}}

\newcommand{\R}{\ensuremath{\mathbb{R}}\xspace}

\newcommand{\CC}{\ensuremath{\mathcal{C}}\xspace}
\newcommand{\N}{\ensuremath{\mathbb{N}}\xspace}
\newcommand{\X}{\ensuremath{X}\xspace}
\newcommand{\G}{\ensuremath{G}\xspace}
\newcommand{\Xs}{\ensuremath{X^*}\xspace}

\newcommand{\Xc}{\ensuremath{[X\mspace{1mu}]}\xspace}
\newcommand{\Xf}{\ensuremath{\langle X\rangle}}
\newcommand{\InI}{\ensuremath{\In(I)}\xspace}
\newcommand{\I}[1]{\ensuremath{\mathcal{I}(#1)}\xspace}
\newcommand{\M}{\ensuremath{M}\xspace}
\newcommand{\Iaw}{\ensuremath{I(A,W)}\xspace}
\newcommand{\pr}{\ensuremath{\prec}\xspace}
\newcommand{\era}[2]{\ensuremath{#1 \raisebox{1.5pt}{\(\scriptscriptstyle\backslash\)} #2}\xspace}

\DeclareMathOperator{\In}{\mathit{In}}

\newcommand{\Grb}{Gr\"obner\xspace}
\newcommand{\pin}[1]{\ensuremath{\pi^{-1}\!#1}\xspace}
\newcommand{\pinp}[1]{\ensuremath{\pi^{-1}\!\left(#1\right)}\xspace}
\newcommand{\si}[1]{\ensuremath{\sigma(#1)}\xspace}
\newcommand{\ha}[1]{\ensuremath{\hat{#1}}\xspace}
\newcommand{\su}[1]{\ensuremath{\underline{#1}}\xspace}
\newcommand{\mi}[1]{\ensuremath{\min(\underline{#1})}\xspace}
\newcommand{\ma}[1]{\ensuremath{\max(\underline{#1})}\xspace}
\newcommand{\ideal}[1]{\ensuremath{\left\langle#1\mspace{1mu}\right\rangle}\xspace}

\renewcommand{\int}[1]{\ensuremath{\iota(#1)}\xspace}
\newcommand{\nint}[1]{\ensuremath{[\iota(#1)]}\xspace}

\newcommand{\true}{\emph{true}\xspace}
\newcommand{\false}{\emph{false}\xspace}
\newcommand{\di}[2]{\ensuremath{#1\!\mid\!#2}}
\newcounter{numex}
\newcommand{\example}{\medskip\noindent
  \textsc{Example \thenumex\stepcounter{numex}:\xspace}}

\title{From monomials to words to graphs}

\author{Cristina G. Fernandes\thanks{Research MCT/CNPq through ProNEx
  Programme (Proc.\ CNPq 664107/1997--4).} \\[-5pt]
  \emph{\small Departamento de Ci\^encia da Computa\c{c}\~ao, Universidade de
  S\~ao Paulo,}\\[-5pt] 
  \emph{\small S\~ao Paulo, SP, Brazil 05508-970}\\[-5pt]
  \texttt{\small cris@ime.usp.br}
\and
Edward L. Green\thanks{Research partially supported by CNPq/NSF
  grant n.\ 910123/99} \\[-5pt]
  \emph{\small Mathematics Department, Virginia Tech University}\\[-5pt] 
  \texttt{\small green@math.vt.edu}
\and
Arnaldo Mandel\thanks{Research partially supported by CNPq/NSF
  grant n.\ 910123/99 and MCT/CNPq through ProNEx
  Programme (Proc.\ CNPq 664107/1997--4).} \\[-5pt]
  \emph{\small Departamento de Ci\^encia da Computa\c{c}\~ao, Universidade de
  S\~ao Paulo,}\\[-5pt] 
  \emph{\small S\~ao Paulo, SP, Brazil 05508-970}\\[-5pt]
  \texttt{\small am@ime.usp.br}}
\date{\today}

\begin{document}
\maketitle

{\makeatletter
\renewcommand{\@makefnmark}{}
\makeatother
\footnote{MSC Classification (2000): 13P10 (Primary), 68R15,
  05C17, 16Z05, 68Q25}}

\begin{abstract}
  Given a finite alphabet \X and an ordering \pr on the letters,
  the map \(\sigma_{\!\pr}\) sends each monomial on \X to the
  word that is the ordered product of the letter powers in the
  monomial.  Motivated by a question on \Grb bases, we
  characterize ideals \(I\) in the free commutative monoid (in
  terms of a generating set) such that the ideal
  \ideal{\sigma_{\!\pr}(I)} generated by \(\sigma_{\!\pr}(I)\) in
  the free monoid is finitely generated.  Whether there exists an
  \pr such that \ideal{\sigma_{\!\pr}(I)} is finitely generated
  turns out to be NP-complete.  The latter problem is closely
  related to the recognition problem for comparability graphs.
\end{abstract}

\section{Introduction}
\label{sec:intro}

An important structural difference between commutative and
noncommutative free monoids is that in a finitely generated free
commutative monoid all ideals are finitely generated (Dickson's
Lemma), while this is not the case for free monoids (for
instance, the ideal generated by \(\{xy^nx\mid n\geq 0\}\)).  We will
consider questions about whether some ideals of the free monoid
are finitely generated.  Those ideals will be described starting
with an ideal in a free commutative monoid.


Let \X be a finite alphabet (of \emph{letters}).  Denote by \Xc
the free commutative monoid on \X and by \Xs the free monoid on
\X; we call the members of \Xc \emph{monomials}, and the members
of \Xs \emph{words}.  When convenient, we assume \(X =
\{x_1,x_2,\ldots,x_n\}\), so a monomial can be written
multiplicatively as \(x_1^{i_1}x_2^{i_2}\cdots x_n^{i_n}\).  We
denote by \(\pi\) the canonical monoid epimorphism
\(\Xs\rightarrow\Xc\).

The most natural relation between ideals of \Xc and \Xs is given
by the canonical map, so the first question is:

\begin{problem}\label{pinfg}
  Given an ideal \(I\) of \Xc, is \pin{(I)} finitely-generated?
\end{problem}



This first question, while quite natural, seems to be of limited
interest, as ideals are not such big players in the structure
theory of monoids.  

The next questions were motivated by the study of noncommutative
presentations of affine algebras and their \Grb bases.  In fact,
the problems studied in this paper were motivated by a question
posed by Bernd Sturmfels on the finite generation of monomial
ideals.  We postpone the discussion until Section \ref{sec:ring},
as the questions can be completely understood within the context
of monoids.

For each ordering \pr of the letters, we define a section
\(\sigma_{\!\pr}\) of \(\pi\) as follows.  We say that a word is
\emph{sorted} if its letters occur in it in increasing order; if
\(x_1\pr x_2\pr\cdots\pr x_n\), such a word can be uniquely
written as \(x_1^{i_1}x_2^{i_2}\cdots x_n^{i_n}\).  Let
\(\sigma_{\!\pr}:\Xc\rightarrow\Xs\) be the function mapping each
monomial \(m\) to the unique sorted word in \pinp{m}.  So,
\(\pi\sigma_{\!\pr}\) is the identity map on \Xc, and we define
the \emph{sorting map} \(S_{\!\pr}\) on \Xs by
\(S_{\!\pr}=\sigma_{\!\pr}\pi\).  The subscript \pr will be
omitted when implicitly understood.

We will mainly be concerned with ideals of form
\(\mathcal{I}_{\!\pr}(I)=\bigl\langle\sigma_{\!\pr}(I)\bigr\rangle\),
that is, the ideal generated by sorted words corresponding to a
commutative ideal.


\begin{problem}\label{IMfg}
  Given an ideal \(I\) of \Xc and an ordering \pr
  of \X, is \(\mathcal{I}_{\!\pr}(I)\) finitely-generated?
\end{problem}


While it is convenient to assume an ordering on the letters so
that one can write monomials, that ordering is not part of
\Xc. So, we also consider

\begin{problem}\label{existscool}
  Given an ideal \(I\) of \Xc, is there an ordering \pr of \X so that
  \(\mathcal{I}_{\!\pr}(I)\) is finitely-generated?
\end{problem}


These problems are very loosely posed, as it is not specified how
each ideal is given.  We are interested in specification by
finite data, so that it makes sense to look for an algorithmic
answer to each of the problems.  We consider three forms of
specifying an ideal of \Xc: by a finite generating set, as the
inverse image of an ideal  under a morphism from \Xc, and as the
initial ideal of a polynomial ideal given by its generators.  The
latter is explained in Section \ref{sec:ring}; our main results
relate to the other two.

Most of this paper will be concerned with a finite set of
monomials and the ideal it generates.  Let us denote by \ideal{M}
the monoid ideal generated by a set \M.  So, \(\ideal{M}=M\Xc\)
if \(M\subseteq\Xc\), and \(\ideal{M}=\Xs M\Xs\) if
\(M\subseteq\Xs\).  From now on, except where explicitly
stated, in each of the problems above, \textit{given an ideal
  \(I\)} is reinterpreted as \textit{given a finite set \M of
  monomials, let \(I=\ideal{M}\)}.  That is the first of the
three forms mentioned above.

It is not clear \textit{a priori} that either problem is even
decidable.  Section~\ref{sec:automata} shows that standard
methods of Automata Theory suffice to decide each instance of the
problems, when the ideal is given by generators.  However, this
is unsatisfying both mathematically and computationally.  From a
purely mathematical viewpoint, the automata decision process is so
far removed from the initial data that one learns very little
about the underlying structure.  From the computational
viewpoint, what goes wrong is the exponential complexity of the
algorithms thus obtained.  We heed the fact that a monomial can
be given as a vector of exponents, so the run-time of algorithms
for the corresponding decision problems should be measured
relative to the bit size of the exponents.  With that in mind,
we summarize the main results:

Problem \ref{pinfg} is solved completely in Section
\ref{sec:pre-image}.  The characterization we obtain yields a
polynomial algorithm.

Problem \ref{IMfg} is the central one here, and also has a
definite answer.  To describe it we need more notation.  We
suppose an ordering of \X is given.  For \(w\) in \Xc, say that a
letter is \emph{extremal} in \(w\) if it is the smallest or the
largest letter with a positive exponent there and say that a
letter is \emph{internal} to \(w\) if it lies strictly between
the extremal letters.  Notice that an internal letter is not
required to occur in \(w\); for instance, using the ordering
implied by the indices, the internal letters of
\(x_2^3x_3x_5^2x_7\) are \(x_3\), \(x_4\), \(x_5\), and \(x_6\).
Also denote by \era{w}{x} the monomial resulting of
evaluating~\(x\) to~\(1\) in the monomial~\(w\).  In particular,
\era{x^m}{x=1}.  A collection of monomials is an \emph{antichain}
if no one divides another; clearly, the (unique) minimal
generating set of an ideal of~\Xc is an antichain.  Dickson's
Lemma, quoted earlier, is equivalent to the statement that every
antichain of monomials is finite.  For convenience, we will
shorten \I{\ideal{M}} to \I{M} when \M is an arbitrary set of
monomials.

\begin{teo}\label{fg}
  Let \M be an antichain in \Xc.  Then, \I{M} is finitely
  generated if and only if, for every \(w\) in \M and \(x\) in
  \X, there exists \(s\) in \M such that \(x\) is extremal in
  \(s\) and \era{s}{x} divides \(w\).
\end{teo}


A proof is found in Section \ref{sec:approach}.
The above result immediately yields a polynomial-time decision
algorithm for Problem \ref{IMfg}.  It also implies that Problem
\ref{existscool} is in NP.

When \M is square free, Problem \ref{existscool} can be decided
in polynomial time. And that is the end of good news.  Problem
\ref{existscool} is shown in Section \ref{sec:what-can-be} to be
NP-complete even when \M consists only of quadratic monomials.
So, while the automata-theoretic algorithms where unacceptable
for Problems~\ref{pinfg} and~\ref{IMfg} because of high exponents
in the data, NP-completeness of Problem \ref{existscool} is not
related to the possibility of writing numbers succinctly.

In Section \ref{sec:comm-ideals-given}, using a slightly
homological flavor, we turn to another way of presenting an ideal
of \Xc: fix a homomorphism from \Xc to another free commutative
monoid, and an ideal \(J\) in the target, and take the pre-image
of that ideal.  The homomorphism can be described by a matrix,
and when~\(J\) is given by its minimal generators the ideal of
\Xc is described as a finite union of integer polyhedra, each one
of them an ideal.  The theory developed for handling generators
is enough to show that Problems \ref{pinfg} and \ref{IMfg} become
coNP-complete with this data, while Problem \ref{existscool} is
shown to be NP-hard.

\section{Connections with \Grb bases}
\label{sec:ring}

Problems \ref{IMfg} and \ref{existscool} stem from a connection
between the commutative and noncommutative \Grb bases theories.
Let \(K\) denote a field, \(K[X]\) the commutative polynomial
ring on the finite set \(X\), and \(K\Xf\) the free associative
algebra on the same set (we adhere to the terminology of the
first section, and talk about \emph{letters} instead of
variables).  The linear extension of the monoid morphism \(\pi\)
is a ring morphism \(K\Xf\rightarrow K[X]\), still denoted by
\(\pi\); its kernel is generated by the commutation relations
\(\CC=\{xy-yx\mid x,y\in\X\}\).  Also, given an ordering of \X,
the linear extension of the maps \(\sigma\) and \(S\) will be
denoted by the same symbols.

Throughout this section, \(I\) is an ideal of \(K[X]\) and
\(J=\pinp{I}\).  It is occasionally useful to lift a commutative
ring presentation \(K[X]/I\) to a noncommutative presentation
\(K\Xf/J\) through \(\pi\).  This has been used in \cite{An, AR,
  GH, PRS} for homological computations.

\begin{pro}\label{genJ}
  Let \(A\subseteq\Xf\) be a set of noncommutative polynomials.
  The following are equivalent:
  \begin{itemize}
  \item [i)] \(\CC\cup A\) generates \(J\).
  \item [ii)] \(\pi(A)\) generates \(I\).
  \item [iii)] For some ordering \pr of \X, \(\CC\cup
    S_\pr(A)\) generates \(J\).
  \item [iv)] For any ordering \pr of \X, \(\CC\cup
    S_\pr(A)\) generates \(J\).
  \end{itemize}
  In particular, for any ordering \pr of \X,
  \(\CC\cup\sigma_{\!\pr}(I)\) generates \(J\).
\end{pro}
\begin{proof}
  Suppose that \(\CC\cup A\) generates \(J\).  Then,
  \(\ideal{\pi(A)} = \ideal{\pi(\CC\cup A)} = \pi\ideal{\CC\cup
    A} = \pi(J) = I\).  Conversely, suppose that \(\pi(A)\)
  generates \(I\).  Since \(\ker\pi\subseteq\ideal{\CC\cup
    A}\subseteq\pinp{I}\) and \(\pi\ideal{\CC\cup A}=I\), it
  follows that \(\ideal{\CC\cup A}=J\).

  The equivalence of conditions (\textit{iii}) and
  (\textit{iv}) to the previous ones follows from the
  observation that \(\pi S_{\!\pr}(A)=\pi(A)\).  The last observation is
  immediate, as \(\pi\sigma_{\!\pr}(I)=I\).
\end{proof}

In many applications, one wants to describe a \Grb basis for
\(J\), preferably related to a \Grb basis for \(I\).  Let us
recall quickly what those bases are (see \cite{Gr,Mo} for an
introduction to the subject).  We will say \emph{term} to mean
either ``word'' or ``monomial'', so we can treat commutative and
noncommutative polynomials simultaneously.

A necessary ingredient for a \Grb bases theory is to fix a
\emph{term order}: a total order on the terms, compatible with
multiplication and with \(1\) as minimum, and with no infinite
descending chains.  In this case, every polynomial has an
\emph{initial term}, the maximum term in its support, and to every
ideal \(I\) one associates its \emph{initial ideal} \InI, the set
of all leading terms of polynomials in \(I\).  The initial ideal
is an ideal of the monoid of terms.  A \emph{\Grb basis} for
\(I\) is a subset \(B\) of \(I\) such that the initial terms of
the members of~\(B\) generate \InI.

Eisenbud, Peeva and Sturmfels \cite{EPS} took the following
approach to relate \Grb bases for commutative polynomial ideals
and their pre-images in the free algebra.  Start with a term
order \pr on \Xc, and define its \emph{lexicographic
  extension}, still denoted here by \pr, to \Xs by:
\(u\pr v\) if \(\pi(u)\pr\pi(v)\) or \(\pi(u)=\pi(v)\) and
\(u\) precedes \(v\) lexicographically, according to the
\pr ordering on \X.

\begin{pro}\label{geniniJ}
  The initial ideal \(\In(J)\) is generated by
  \[\{xy\mid x,y\in X, \; x\succ y\}\cup\si{\InI}.\]
\end{pro}
\begin{proof}
  Since the noncommutative order extends the commutative one, if
  \(p\in K[X]\), the initial term of \si{p} is \si{m}, where
  \(m\) is the initial term of~\(p\).  Also, from the
  lexicography, if \(x\succ y\) are letters, \(xy\) is the
  initial term of \(xy-yx\).  The result now follows from
  Proposition \ref{genJ}.
\end{proof}

It follows from Dickson's Theorem that every ideal of \Xc is
finitely generated, hence every ideal of \(K[X]\) has a finite
\Grb basis.  In contrast, not every ideal of \Xs is finitely
generated, so it is generally interesting to detect whether a
given ideal of \(K\Xf\) has a finite \Grb basis.  The preceding
proposition implies that \(J\) has a finite \Grb basis with
respect to \(\succ\) if and only if \I{\InI} is finitely
generated.  This gives rise to Problem \ref{IMfg}.

A word in an ideal of \Xs is a minimal generator of that ideal if
and only if the words obtained by erasing either the first or the
last letter are not in the ideal.  Combining it with Proposition
\ref{genJ}, we get the next result; it is Theorem 2.1 of
\cite{EPS}, stripped of the ring theoretic context (which is
handled by Proposition~\ref{geniniJ}):

\begin{teo}
  If \M is a multiplicative antichain of monomials, then the
  minimal generating set of \I{M} is
  \vspace{-2ex}
    \begin{multline*}
      \{\si{mu} \mid m\in M, \text{\(u\in\Xc\) is generated by
        letters internal to \(m\),}\\
      \text{ and is such that, for each letter \(x\) extremal to
        \(m\), } mux^{-1}\not\in\ideal{M}\}.
    \end{multline*}
\end{teo}


There seems to be no immediate characterization from the above
for when~\(J\) has a finite \Grb basis.  A sufficient condition
is provided in~\cite{EPS}, and then the question is finessed: it
is shown that, if \(K\) is infinite, then, for any \(I\) and
\pr, \(J\) will have a finite \Grb basis after a generic
change of variables,

This is of limited computational use if high degree polynomials
are being handled, since the supports can grow explosively after
generic changes.  So, it is still of some interest to detect
whether \pinp{I} has a finite \Grb basis when \(I\) is still
expressed in the given coordinates.

If \(I\) is generated by a set \M of monomials (this is
called a \emph{monomial ideal} of \(K[X]\)), then \InI =
\ideal{M}, irrespective of the term order \pr.  It follows
from Proposition \ref{geniniJ} that:

\begin{pro}\label{anyJ}
  Let \(M\subseteq \Xc\) and let \(I\) be the monomial ideal it
  generates.  Then \(J\) has a finite \Grb basis with respect to
  the lexicographic extension of a term order of \Xc if and only
  if, for the same ordering of \X, \I{M} is finitely generated.
\end{pro}

This gives rise to Problem \ref{existscool}.

A similar looking question is, in the notation above, what
conditions must \M satisfy, so that \(J\) has a finite \Grb basis
with respect to the lexicographic extension of any term order of
\Xc?  Proposition \ref{anyJ} translates it to a problem of
monomials and words, and the answer is in section \ref{sec:p2},
Theorem \ref{totallycool}.

As a final note, we point out that within this context another
way of specifying an ideal of \Xc is relevant to Problems
\ref{IMfg} and \ref{existscool}.  Namely, suppose a term order is
given on \Xc; given a finite set \M of polynomials, consider the
initial ideal \(I\) of the polynomial ideal generated by \M.
That is the actual motivation for those problems, after all!
The usual process of going from \M to \(I\) is Buchberger's
algorithm and its variants.  These are all of high complexity,
so the question remains whether Problem \ref{IMfg} can be solved
efficiently from this data.

\section{Using automata}
\label{sec:automata}

It is not clear from the outset that either of the problems
mentioned in the introduction is decidable.  This can be shown to
be the case by means of the traditional machinery of automata
theory (we follow the notation and terminology of \cite{LP}).  We
do it here, mostly for completeness and to underline some of the
complexity issues.  This section can be skipped, with no loss in
understanding of the remaining text.

Suppose \(J\) is an ideal of the free monoid \Xs; then its unique
minimal generating set is \(T=J\backslash(XJ\cup JX)\).  Hence,
if \(J\) is a regular language, so is~\(T\).  The problem of
deciding whether \(T\) is finite, given a regular expression for
\(J\), can be whimsically, although not very accurately, related
to two well-known Unix utility programs: given a pattern for a
\texttt{grep} search, decide whether the same search can be made
by \texttt{fgrep}.  

Problems \ref{pinfg} and \ref{IMfg} refer to ideals that are
regular languages.  Consider Problem \ref{pinfg}. It is a
well-known (although nonconstructive) consequence of Higman's
Theorem \cite{Hi} that, for every subset \(A\) of \Xc, \pinp{A} is
regular.  For \pin{\ideal{M}}, one can construct a deterministic
automaton directly: have \(n\) parallel counters, one for each
letter and counting up to its maximum degree in \M.  Each state
of the automaton corresponds to an \(n\)-tuple of values for the
counters, and processing a word \(x\) leads to a state whose
counters correspond to the exponent vector of \(\pi(x)\).  The
final states are those that show that \(\pi(x)\in\ideal{M}\), and
they can be colluded into a single state that is never left after
being reached.

For Problem \ref{IMfg}, we can write a simple regular
expression for \pin{\ideal{M}}.
If \(w=x_1^{i_1}x_2^{i_2}\cdots x_{n-1}^{i_{n-1}}
x_n^{i_n}\in\Xc\), then \(\I{w}=\Xs
x_1^{i_1}x_2^{i_2}x_2^*\cdots
x_{n-1}^{i_{n-1}}x_{n-1}^*x_n^{i_n}\Xs\), and \(\I{M}=\cup_{w\in
  M} \I{w}\), so \I{M} is a regular language.  A
deterministic automaton recognizing \I{M} can also be constructed
using parallel bounded counters, although the description would
be more complicated than the previous one.

In both cases, there is only one final state, from which no
transition leaves.  This makes it easy to construct a
deterministic automaton for the minimal generating set of each
ideal, with direct products of three very similar automata.  Then
finiteness can be easily checked by a graph search.  This
approach shows now that Problem \ref{IMfg} is decidable.

From a complexity viewpoint, this does not work.  Even though we
have scrupulously avoided using nondeterministic automata, there
remains a source of exponential complexity:  in either case, the
automaton described for each ideal has a number of states that is
roughly the product \(N\) of the maximum degrees of letters in
\M.  This is too large;  since a monomial can be represented as a
vector of exponents, a reasonable encoding for \M would have only
\(O(n\log N)\) bits, where \(n=|X|\).  So, the automaton for the
minimal generating set has exponentially many states, and the
graph search is linear in the number of states.

The solutions we present in the following sections could perhaps
be retro\-fit\-ted into an automata-theoretical framework.
Actually, thinking of automata helped in the discovery of those
results: the Pumping Lemma (see \cite[Sec. 2.4]{LP}) was a
starting point.

There is a rich literature on ideals of the free monoid, from the
viewpoint of language theory, with a twist.  Instead of
concentrating on the ideal, the focus is on its complement.
The complement of an ideal is said to be  \emph{factorial
  languages}, and the minimal generators of the ideal appear as
\emph{forbidden subwords} in this context.  There are many
algorithms for problems involving factorial languages (see \cite{BMRS,
CFR}), but they usually take a deterministic automaton like the
large ones we described as input, so they are of no use here.

\section{When \I{M} is finitely generated}
\label{sec:approach}

We assume a fixed ordering \pr on the alphabet \X. The
\emph{support} of a monomial~\(w\), denoted \su{w}, is the set of
letters with nonzero exponent in \(w\).  So, \mi{w} and \ma{w}
are the \emph{extremal} letters of \(w\), while the letters
\(x\) such that \(\mi{w} \pr x \pr \ma{w}\) are \emph{internal} to
\(w\).  We will use the notation \di{u}{v} meaning \emph{\(u\)
  divides \(v\)}, both in \Xc and \Xs.  So, for monomials
\(u=x_1^{i_1}\ldots x_n^{i_n}\), \(v=x_1^{j_1}\ldots x_n^{j_n}\),
\di{u}{v} means that \(i_1\leq j_1,\ldots,i_n\leq j_n\); for
words \(u\), \(v\), \di{u}{v} means that there exist words
\(w\), \(z\) such that \(v=wuz\).

\example Take \(X=\{a,b,c,\dots\}\) with the usual ordering.
Then, for \(w=b^2df\), \(\mi{w}=b\), \(\ma{w}=f\), and the
internal letters are \(c\), \(d\), and \(e\).

\begin{lem}\label{divisor}
  If \(u\) and \(v\) are sorted words such that \di{u}{v}, then
  for any \(x\) in \(X\), either \di{u}{S(vx)} or
  \di{S(ux)}{S(vx)}.
\end{lem}
\begin{proof}
  The first case occurs if \(x\leq\mi{u}\) or \(x\geq\ma{u}\),
  and the second case occurs if \(\mi{u}\leq x \leq \ma{u}\).
\end{proof}

\begin{lem}\label{generator}
  Let \(T\) be a set of sorted words.  If, for every \(t\) in
  \(T\) and \(x\) in \X, \(S(tx)\) has a factor in \(T\), then
  there exists an ideal \(I\) of \Xc such that \(T\subseteq
  \si{I}\) and \(\ideal{T} = \I{I}\).
\end{lem}
\begin{proof}
  Let \(M=\pi(T)\) and \(I=\ideal{M}\); clearly \(T\subseteq
  \si{I}\), and we will show that \(\ideal{T} = \I{I}\).
  
  Since \(\ideal{T} \subseteq \I{I}\) is clear, it is enough to
  show that \(\si{I} \subseteq \ideal{T}\).  Since any monomial
  in \(I\) can be written as \(su\), with \(s\in M\),
  \(u\in\Xc\), we will show that \(\si{su}\) has a factor in
  \(T\) (so, it is in \ideal{T}) by induction in the total degree
  of \(u\).
  
  There is nothing to prove if the degree is zero.  So, we can
  write \(u=vx\), with \(v\in\Xc\), \(x\in X\).  By the induction
  hypothesis, \di{t}{\si{sv}} for some \(t\in T\).  By
  hypothesis, \(S(tx)\) has a factor \(y\in T\).  Noticing that
  \(S(\si{sv}x)=\si{svx}=\si{su}\), Lemma \ref{divisor} implies
  that either \di{t}{\si{su}}, or \di{S(tx)}{\si{su}}, in which
  case, \di{y}{\si{su}}.  In either case, the result follows.
\end{proof}

From this, it easily follows:
\begin{cor}\label{crit}
  Let \M be a set of monomials.  A sufficient condition for a set
  \(T\) of words to be such that \(\ideal{T} = \I{M}\) is that
  \(\si{M}\subseteq T\subseteq \sigma\!\ideal{M}\) and for every
  \(t\) in \(T\) and \(x\) in \X, \(S(tx)\) has a factor in \(T\).
\end{cor}

For the proof, we restate Theorem \ref{fg} with some additional
precision.  First, we recall and introduce some notation.

If \(w\) is a monomial and \(x\in X\), \era{w}{x} denotes the monomial
obtained from~\(w\) by erasing the occurrences of \(x\).  We denote by
\int{w} the set of internal letters of \(w\), and by \(\partial_xw\)
the degree of \(x\) in \(w\).  Given a set \M of monomials, let
\(r_x(M)\) denote the maximum degree \(x\) occurs with as an extremal
letter in \M.  To avoid misunderstandings, \nint{w} is simply the
submonoid of \Xc generated by~\int{w}.

\example Continuing the earlier example, \(\era{w}{b}=df\) and
\(\era{w}{f}=b^2d\).  If \(M=\{c^3,a^2c^5f^2, cf^3g,a^2b^2c^2\}\),
then \(r_c(M)=3\) and \(r_f(M)=2\).

\begin{teo}\label{FG}
  Let \M be an antichain in \Xc.  The following are equivalent:
  \begin{itemize}
  \item[i)] \I{M} is finitely generated.
  \item[ii)] For every \(w\) in \M and \(x\) in \int{w}, there
    exists  \(s\) in \M such that \(x\) is extremal in \(s\), and
    \era{s}{x} divides \(w\).
  \item[iii)] \I{M} is generated by
    \(\displaystyle\sigma\left(\smash[b]{\bigcup_{w\in M}} \{u\in
      w\nint{w}\mid \forall x\in \int{w},\ 
      \partial_xu<r_x(M)\}\right).\)
  \end{itemize}
\end{teo}

  
\begin{proof}
  Let \(T\) be the minimal generating set of \I{M}, and suppose
  it is finite.  We shall prove condition (\textit{ii}). Note
  that \(\si{M}\subseteq T\).  Indeed, if \(w\in M, \si{w}\) has
  a factor in \(T\), and this has the form \(\si{su}\), for some
  \(s\in M\) and \(u\in\Xc\).  So, \di{\si{su}}{\si{w}}, and this
  clearly implies \di{su}{w}.  So, \di{s}{w} and, since \M is an
  antichain and \(s\in M\), it follows that \(s=w, u=1\), hence
  \(\si{w}\in T\).
  
  Let \(w\in M\) and let \(x\) be an internal letter to \(w\).
  We can uniquely write \(\si{w}=ux^rv\), with \(u, v\) sorted,
  \(\ma{u} \pr x \pr\mi{v}\).  Hence, all words \(ux^nv\), with
  \(n\geq r\), are in \I{M}, so each has a factor in \(T\); it
  follows that some \(t\in T\) is a factor of infinitely many
  such words.  If \(x\) is not in the support of \(t\), it must
  happen that \(t\) is a factor of \(u\) or of \(v\), hence a
  proper factor of \(\si{w}\), but \(\si{w}\in T\), so this
  cannot occur.  Therefore, \(x\) is in the support of \(t\), and
  necessarily is extremal.  Without loss of generality, let us
  assume that \(x=\mi{t}\); so \(t=x^kz\) with \(x\pr\mi{z}\), and
  note that \(z\) is a factor of \(\si{w}\), so \di{\pi(z)}{w}.
  Now, \(t=\si{sy}\) for some \(s\in M\) and \(y\in \Xc\).  If
  \(x\notin \su{s}\), we would have \di{s}{\pi(z)}, hence
  \di{s}{w}, a contradiction.  Hence \(x=\mi{s}\), so \era{s}{x}
  is a factor of \(w\).
  
  
  Now, suppose that condition (\textit{ii}) holds and let us
  prove (\textit{iii}).  Call \(T\) the generating set in that
  statement.  Clearly \(\si{M}\subseteq T\subseteq
  \sigma\!\ideal{M}\).  Now, let \(t\in T\) and \(x\in X\), and let
  us find a factor of \(S(tx)\) in \(T\) as required by Corollary
  \ref{crit}.  Write \(t=\si{wz}\), with \(w\in M\) and \(z\in
  \nint{w}\).
  
  If \(x\notin\int{w}=\int{\pi(t)}\), it follows immediately that
  \di{t}{S(tx)}.  There remains the case where \(x\in\int{w}\)
  and \(S(tx)\not\in T\) (since the case \(S(tx)\in T\) is
  trivial).
  
  Now, \(\si{wz}\in T\), but \(\si{wzx}\notin T\), hence
  \(\partial_xwz=r_x(M)-1\).  We can write uniquely
  \(wz=ux^{r_x(M)-1}v\), with \(\ma{u} \pr x \pr\mi{v}\).  By
  hypothesis, there exists \(s\in M\), with \(x\) extremal in
  \(s\) (without loss, \(x=\ma{s}\)) such that \(\era{s}{x}\)
  divides \(w\).  Since \(\partial_xs\leq r_x(M)\), \di{s}{wzx},
  and by maximality of \(x\) in~\su{s}, \di{s}{ux^{r_x(M)}}.  Let
  \(p\) result from raising the degree of each internal letter of
  \(s\) to its exponent in \(u\).  Then, \(p\in s[\int{s}]\) and
  its internal letters have small degree, so \(\si{p}\in T\), and
  it is the factor of \(S(tx)\) we sought after.
  
  Clearly (\textit{iii}) implies (\textit{i}), and the
  theorem is proved.
\end{proof}

\example Consider the set \(M=\{ab^2c,a^3b\}\).  The ordering
\(a\pr b\pr c\) does not satisfy the conditions above, since
\(b\) is internal to \(ab^2c\), and extremal only in \(a^3b\),
but \(\era{a^3b}{b}=a^3\) does not divide \(ab^2c\); indeed,
\(\sigma\!\ideal{M}\) comprises all words \(a^ib^{j+2}c^k\) and
\(a^{i+3}b^j\) with \(i,j,k\geq 0\), and the set \(\{ab^{j+2}c
\mid j\geq 0\}\) cannot be generated as multiples of finitely
many of those.  Similarly, \(a\pr c\pr b\) fails, since \(c\) is
internal to \(a^3b\), and extremal in none.  However, \(b\pr a\pr
c\) is good: \(a\) is internal to \(b^2ac\) only, and
\era{ba^3}{a=b} divides \(b^2ac\).  In this case, \I{M} is
generated by \(\{b^2ac,b^2a^2c,ba^3\}\).

Condition (\textit{iii}) above is a fairly precise description
of the minimal generating set of \I{M}.  One gets a quick and
dirty estimate for its size by forgetting most parameters:

\begin{cor}
  Let \M be a finite set of monomials in \Xc.  Then, if \I{M} is
  finitely generated, it can be generated by a set of at most
  \[|M'|\prod_{x\in\int{M}}r_x(M)+|M|-|M'|\]
  elements, where
  \(M'=\{w\in M\mid \int{w}\neq\emptyset\}\) and \(\int{M}=\bigcup_{w\in
    M}\int{w}\).
\end{cor}


Even though it is a rough estimate, the result above is best
possible.  To see this, suppose \(X = \{x_1,x_2,\ldots,x_n\}\),
ordered according to the indices.  Choose positive integers
\(m,r_2,r_3,\ldots,r_{n-1}\), and let
\(M=\{x_1^{i+1}x_n^{m-i}\mid 0\leq i <m\}\cup\{x_i^{r_i}\mid
i=2,\ldots,n-1\}\).  Then, the minimal generating set for \I{M}
is precisely that described in Theorem \ref{FG} (\textit{iii})
and has size \(m\prod r_i+n-2\).

If we are given \M as a collection of integer vectors,
divisibility is just componentwise comparison, so it can be
tested rapidly.  A naive check of the condition on Theorem
\ref{fg} would need at most \(|M|^2|X|\) such comparisons, so
Problem \ref{IMfg} can be solved by a polynomial-time algorithm.


We end this section with some constructions that will be useful
later and some unexpected consequences of the theorem.

Given a collection \M of monomials and an integer \(k\), \(M_k\)
will denote the subset of \M consisting of monomials whose
support has size at most \(k\).

\begin{pro}\label{level}
  If an antichain \M of monomials is such that \I{M} is finitely
  generated, so is \I{M_k} for each integer \(k\).
\end{pro}
\begin{proof}
  Let us show that \(M_k\) satisfies condition (\textit{ii}) of
  Theorem \ref{FG}.  If \(w\in M_k\), and \(x\) is internal to
  \(w\), we know, since \I{M} is finitely generated, that for
  some monomial \(u\in M\), \(x\) is extremal in \(u\) and
  \era{u}{x} is a factor of \(w\).  Since \(x\) is internal to
  \(w\) and extremal in \(u\), the support of \era{u}{x} is a
  proper subset of the support of \(w\), so \(|\su{u}|\leq k\);
  that is, \(u\in M_k\).
\end{proof}

\begin{pro}\label{squash}
  Let \(w\leftrightarrow\ha{w}\) be a bijection between sets \M
  and \ha{M} of monomials, such that:
  \begin{itemize}
  \item[i)] For every \(w\) in \M, \(w\) and \ha{w} have the same
    extremal letters.
  \item[ii)] For every \(u,v\) in \M and \(x\) in \X, if
    \(\partial_xu\leq\partial_xv\), then
    \(\partial_x\ha{u}\leq\partial_x\ha{v}\).
  \end{itemize}
  If \ha{M} is an antichain, and \I{M} is finitely
  generated, then so is \I{\ha{M}}.
\end{pro}
\begin{proof}
  First notice that, from condition (\textit{ii}), \di{u}{v}
  implies \di{\ha{u}}{\ha{v}}, so \M is an antichain.  Thus, it
  must satisfy Theorem \ref{FG} (\textit{ii}).  Let \(\ha{w}\in
  \ha{M}\), and \(x\in X\) be internal to \ha{w}.  From condition
  (\textit{i}), \(x\) is internal to \(w\), so there exists
  \(u\in M\) such that \(x\) is extremal in \(u\) and
  \di{\era{u}{x}}{w}.  Again from (\textit{i}), \(x\) is
  extremal in \ha{u}, and from condition (\textit{ii}),
  \di{\era{\ha{u}}{x}}{\ha{w}}.
\end{proof}

\section{Cool orders}
\label{sec:p2}

Given an antichain \M of monomials in \Xc, we say that an
ordering of \X is \emph{cool} for \M if, for every \(w\) in \M and
letter \(x\) internal to \(w\), there exists \(s\) in \M such
that \(x\) is extremal in \(s\) and \era{s}{x} divides \(w\).

As we will see in the next section, no good algorithm is
forthcoming to decide whether a cool ordering exists.  Before
giving substance to this, we try to get a better understanding of
such orderings.  We begin with an immediate consequence of
Propositions \ref{level} and \ref{squash}:

\begin{pro}\label{coolmon}
  Let \M be an antichain of monomials.  Then, any cool ordering
  for \M is also a cool ordering for:
  \begin{enumerate}
  \item \(M_k\), for any integer \(k\).
  \item Any \ha{M}, obtained from \M by changing each \(w\) to
    \ha{w}, so that:
    \begin{enumerate}
    \item \(w\) and \ha{w} have the same support, and
    \item  For every \(u,v\) in \M and \(x\) in \X, if
    \(\partial_xu\leq\partial_xv\), then
    \(\partial_x\ha{u}\leq\partial_x\ha{v}\).
    \end{enumerate}
  \end{enumerate}
\end{pro}

This gives some necessary conditions for existence of cool
orderings.  We also get a kind of equivalence between sets of
monomials:

\begin{pro}\label{strict}
   Let \(w\leftrightarrow\ha{w}\) be a bijection between set \M
  and \ha{M} of monomials, such that:
  \begin{itemize}
  \item[i)] For every \(w\) in \M, \(\su{w}=\su{\ha{w}}\), and
  \item[ii)] For every \(u,v\) in \M and \(x\) in \X, 
    \(\partial_xu<\partial_xv\) if and only if
    \(\partial_x\ha{u}<\partial_x\ha{v}\).
  \end{itemize}
  Then, an ordering of \X is cool for \M if and only if it is so
  for \ha{M}.
\end{pro}
\begin{proof}
  It is easily checked that \di{u}{v} if and only if
  \di{\ha{u}}{\ha{v}}.  So, the bijection maps minimal
  monomials to minimal monomials, antichains to antichains, and
  so on.  It is just a matter of applying part 2 of Proposition
  \ref{coolmon} in both directions.
\end{proof}

From now to the end of the article, an ordering of the letters
will not be given at the outset, and the following concept will
be useful for the search of cool orderings.  A monomial \(w\) is
said to \emph{help} a monomial \(m\) \emph{with} the letter~\(x\)
if \(x\in\su{w}\), \di{\era{w}{x}}{m} and
\(\su{\era{w}{x}}\subsetneq\su{\era{m}{x}}\).

If \pr is a cool ordering for a set \M and, for some
\(Y\subset X\), every monomial in \M has support included in
\(Y\) or disjoint from \(Y\), then \pr is a cool ordering for
those monomials with support included in \(Y\).  This can be
extended to the following easily verified fact:

\begin{pro}
  A cool ordering for \M is also cool for any \(N\subseteq M\) such
  that all members on \M that help some member of \(N\) (with
  some letter) are in \(N\).
\end{pro}

Next section will consider the problem of finding a cool
ordering, given~\M. We close the section considering a question
that is a sort of opposite of that: how must \M look like if
\emph{every} ordering of \X is cool?  The answer is surprisingly
simple.

\begin{teo}\label{totallycool}
  Let \M be an antichain of monomials over \X.  Then, every
  ordering of \X is cool for \M if and only if, for every \(m\)
  in \M and \(x\) in \X such that \(|\su{\era{m}{x}}|\geq 2\),
  there exists a \(u\) in \(M_2\) such that \di{\era{u}{x}}{m}.
\end{teo}
\begin{proof}
  Suppose that every ordering of \X is cool for \M.  Let \(m\) and
  \(x\) be given as in the statement.  Then, there exists an ordering
  \pr on \X such that~\(x\) is internal to \(m\).  Since \pr is cool,
  there is a \(u\) in \M that helps~\(m\) with~\(x\).  Choose \(u\)
  with minimal support, and let us show that \su{u} has size at most
  \(2\).  If \(|\su{u}|>2\), there is another ordering of \X that
  makes \(x\) internal to \(u\).  Again by Theorem \ref{FG}, there
  exists a \(v\) in \M that helps \(u\) with \(x\); clearly, \(v\)
  also helps~\(m\) with \(x\).  Since \(x\in \su{u}\cap\su{v}\), it
  follows that \(\su{v}\subsetneq\su{u}\), so we have contradicted the
  minimality of \su{u}.
  
  Conversely, suppose the divisibility condition holds, and
  consider an arbitrary ordering of \X.  Pick an \(m\in M\) and
  let \(x\in\int{m}\).  Choose \(u\) in \(M_2\) such that
  \di{\era{u}{x}}{m}; clearly \(x\) is extremal in \(u\).  It
  follows that the ordering is cool.
\end{proof}

\section{Finding cool orders is hard}
\label{sec:what-can-be}

Monoids generated by square-free monomials appear frequently in
algebraic combinatorics (related to Stanley-Reisner rings of
simplicial complexes), and have been studied in the current
context by Peeva and Sturmfels, together with Eisenbud \cite{EPS}
and Reiner \cite{PRS}.  Propositions
\ref{sqfree}~and~\ref{compara} tell the same as \cite[Prop.
3.2]{EPS} and \cite[Lemma 3.1]{PRS}, although the different
jargon may obscure this.  After that we move to another
direction.
 
\begin{pro}
  Suppose that \M is an antichain and the degree of the letter
  \(x\) in \(w\in M\) is the largest degree it has in all
  monomials in \M.  Then, in any cool ordering for \M, \(x\)
  cannot be internal to \(w\).
\end{pro}
\begin{proof}
  If there is a cool ordering for \M where \(x\) is internal to
  \(w\), there exists a \(t\in M\) such that \di{\era{t}{x}}{w}.
  But since \(\partial_xw\geq\partial_xt\), it follows that
  \di{t}{w}, a contradiction, as \M is an antichain.
\end{proof}

The following is an immediate corollary:

\begin{pro}\label{sqfree}
  If \M consists only of square-free monomials and affords a
  cool ordering, then its monomials have total degree at most 2.
\end{pro}

Degree 1 monomials are trivially handled here, so the square-free
sets~\M of interest consist only of quadratic monomials.
Polynomial ideals whose initial ideals are generated by quadratic
monomials (mostly square-free) were extensively studied in
\cite{PRS}.  

We leave now the square-free condition, and consider the case
when \M consists exclusively of quadratic monomials, that is, we
allow monomials of form \(x^2\).  This seemingly trivial
extension has deep consequences:

\begin{pro}\label{qnpc}
  The problem of deciding, given a set of quadratic monomials \M,
  whether there exists a cool ordering for \M is:
  \begin{itemize}
  \item [a) ] Solvable in polynomial time, if \M is square-free.
  \item [b) ] NP-complete, in general.
  \end{itemize}
\end{pro}
\begin{proof}
  Part (a) follows from Proposition \ref{compara} and part (b)
  from Proposition~\ref{npc}.
\end{proof}

With quadratic monomials, irrespective of the order of the
letters, each letter is extremal in each monomial it occurs, so,
an ordering \pr on \X is cool for \M if and only if whenever
\(x\pr y\pr z\) and \(xz\) is in \M, then at least one of \(y^2, xy, yz\)
is in \M.

At this point, it becomes convenient to encode the data and the
problem by means of graphs, and it turns out to be convenient to
use the complement of what comes naturally.  The graph \(G(M)\)
will have the letters as vertices, \(xy\) is an edge if
\(xy\) is not in \M.  Let \(T_M\) denote the set of letters whose
square is \emph{not} in \M.

An orientation of a graph is said to be \emph{transitive} at a
vertex \(y\) if, whenever oriented edges \(x\rightarrow y\) and
\(y\rightarrow z\) exist, then the edge \(x\rightarrow z\) must
also exist.  A graph is a \emph{comparability graph} if it admits
an orientation that is transitive at all its vertices; such an
orientation is always acyclic.  Comparability graphs have been
widely studied, and have efficient recognition algorithms
\cite{Ga} (or \cite {MP}), \cite{Go}, \cite{RAR}. 

\begin{pro}\label{compara}
  A set \M of quadratic monomials admits a cool order if and only
  if \(G(M)\) admits an acyclic orientation that is transitive at
  all vertices of \(T_M\).  In particular, if \M is square-free,
  it admits a cool order if and only if \(G(M)\) is a
  comparability graph.
\end{pro}
\begin{proof}
  Suppose \M has a cool ordering.  Direct all edges of \(G(M)\)
  from the smallest to the largest vertex.  This orientation is
  trivially acyclic.  If \(y\in T_M\), and edges \(x\rightarrow
  y\) and \(y\rightarrow z\) exist, then the monomials \(y^2\), \(xy\)
  and \(yz\) are not in \M.  By coolness, \(xz\not\in M\), so the
  edge \(xz\) is in \(G(M)\), and is correctly oriented.
  
  Conversely, suppose that \(G(M)\) admits an acyclic orientation
  that is transitive at all vertices of \(T_M\).  With a
  ``topological sort'' order its vertices so that all directed
  edges point from the smaller to the bigger end.  One readily
  verifies that this ordering is cool for \M.
  
  When \M is square-free, \(T_M\) comprises all vertices, so an
  acyclic orientation that is transitive at all vertices of
  \(T_M\) says that \(G(M)\) is a comparability graph.
\end{proof}

We refer the reader to the already classic text \cite{GJ} as a
general reference for NP-completeness, good algorithms and
satisfiability.  Since good algorithms for recognition of
comparability graphs are known, one would expect that testing the
condition of Proposition \ref{compara} would also be feasible.

\gasset{Nh=2, Nw=2, Nadjust=n, ExtNL=y, NLangle=90, NLdist=1.2, AHnb=0}
\begin{pro}\label{npc}
  The problem: 
  \begin{quote}
    \emph{given a graph \G and a set \(T\subseteq VG\), is there
      an acyclic orientation of \G that is transitive at all
      vertices in \(T\)?}
  \end{quote}
  is NP-complete.
\end{pro}
\begin{proof}
  Let us shorten ``orientation transitive at \(T\)'' to
  \emph{\(T\)-orientation}.  The proof will be by a reduction
  from \textsc{not-all-equal-3sat} \cite{Sc}.  The basic gadget is
  the graph in Figure \ref{fig:1}, where the vertices in \(T\)
  are black (and labeled \(a\), \(\bar{a}\), \(c\)).

  \begin{figure}[h]
    \centering
    \begin{picture}(70,25)(10,0)
      \node(s)(10,21){\(s\)}
      \node[Nfill=y](a)(30,21){\(a\)}
      \node[Nfill=y](ab)(50,21){\(\bar{a}\)}
      \node(t)(70,21){\(t\)}
      \node[NLangle=180](l)(30,1){\(l\)}
      \node[NLangle=0](r)(50,1){\(r\)}
      \node[Nfill=y](c)(40,11){\(c\)}
      \drawedge(s,a){}
      \drawedge(ab,a){}
      \drawedge(a,l){}
      \drawedge(r,ab){}
      \drawedge(c,ab){}
      \drawedge(c,a){}
      \drawedge(c,r){}
      \drawedge(c,l){}
      \drawedge(t,ab){}
      \drawedge(l,r){}
    \end{picture}
    \caption{The top hat}
    \label{fig:1}
  \end{figure}
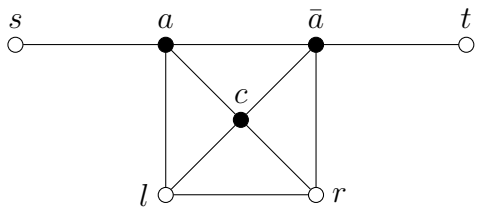
  
  \textsc{Fact 1:} \textit{There is only one \(T\)-orientation of
    this graph where the edge \(a\bar{a}\) is oriented
    \(a\rightarrow\bar{a}\). In that orientation, \(a\) is a
    source, \(\bar{a}\) is a sink, and the bottom edge is
    directed from \(r\) to \(l\).}

  To see this, notice that since the edge \(s\bar{a}\) does not
  exist, \(sa\) must be oriented as \(a\rightarrow s\), because
  of transitivity at \(a\).  By a similar argument we check that
  all edges with an end in \(T\) can have only one orientation.
  Finally, since the orientations \(r\rightarrow c\) and
  \(c\rightarrow l\) are forced, \(r\rightarrow l\) is forced by
  transitivity at \(c\).

  Now we construct the main gadget by gluing three copies of the
  top hat, identifying cyclically each \(t\) with the next \(s\)
  and each \(r\) with the next \(l\).  The result is in Figure
  \ref{fig:2}, where only important vertices are labeled.

  \begin{figure}[h]
    \centering
    \begin{picture}(75,60)(10,0)
      \node(s1)(22.5,42){}
      \node[Nfill=y](a1)(30,55){\(a_1\)}
      \node[Nfill=y](a1b)(60,55){\(\bar{a}_1\)}
      \node(s2)(67.5,42){}
      \node(s3)(45,1){}
      \node[Nfill=y,NLangle=0](a2)(75,29){\(a_2\)}
      \node[Nfill=y,NLangle=0](a2b)(60,1){\(\bar{a}_2\)}
      \node[Nfill=y,NLangle=180](a3)(30,1){\(a_3\)}
      \node[Nfill=y,NLangle=180](a3b)(15,29){\(\bar{a}_3\)}
      \node(l1)(30,37.3){}
      \node(l2)(60,37.3){}
      \node(l3)(45,8.2){}
      \node[Nfill=y](c1)(45,46.2){}
      \node[Nfill=y,NLangle=0](c2)(60,18.5){}
      \node[Nfill=y,NLangle=180](c3)(30,18.5){}
      \drawedge(s1,a1){} \drawedge(a1b,s2){}
      \drawedge(s2,a2){} \drawedge(a2b,s3){}
      \drawedge(s3,a3){} \drawedge(a3b,s1){}
      \drawedge(a1,a1b){}
      \drawedge(a2,a2b){}
      \drawedge(a3,a3b){}
      \drawedge(a1,l1){} \drawedge(a1b,l2){}
      \drawedge(a2,l2){} \drawedge(a2b,l3){}
      \drawedge(a3,l3){} \drawedge(a3b,l1){}
      \drawedge(c1,a1){} \drawedge(c1,a1b){} \drawedge(c1,l1){} \drawedge(c1,l2){} 
      \drawedge(c2,a2){} \drawedge(c2,a2b){} \drawedge(c2,l2){} \drawedge(c2,l3){} 
      \drawedge(c3,a3){} \drawedge(c3,a3b){} \drawedge(c3,l3){}
      \drawedge(c3,l1){}
      \drawedge(l1,l2){} \drawedge(l2,l3){} \drawedge(l3,l1){} 
    \end{picture}
    \caption{The gadget}
    \label{fig:2}
  \end{figure}
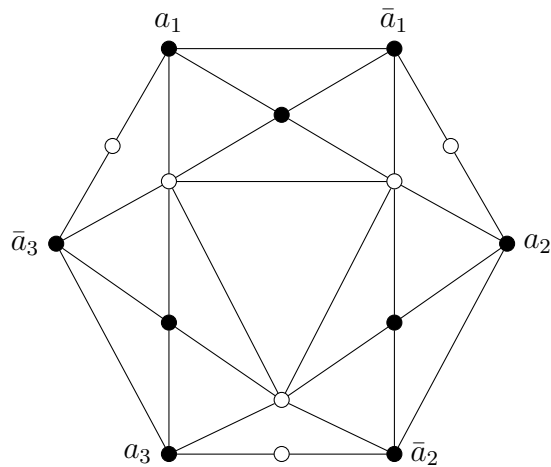
  
  \textsc{Fact 2:} \emph{Consider an orientation of the edges
    \(a_1\bar{a}_1\), \(a_2\bar{a}_2\) and \(a_3\bar{a}_3\).  It
    extends to an acyclic \(T\)-orientation of the gadget if and
    only if they are not all directed the same way along the
    external cycle.}
  
  Indeed,  by looking at the top hats we see that any
  orientation of these edges extends uniquely to a
  \(T\)-orientation of the gadget.  If they are all oriented the
  same way, the inner triangle becomes a directed cycle.
  Conversely, if they are not all the same way (by symmetry,
  there is only one case to check), the orientation of the gadget
  is acyclic.

  Now we proceed to the reduction.  A typical instance of
  \textsc{not-all-equal-3sat} consists of a set \X of
  variables and a set \(C\) of clauses over \X, where each
  clause has three literals, each of form \(x\) or \(\bar{x}\),
  for some \(x\in X\).  The question is whether there exists a
  truth assignment to \X, so that for each clause  one
  literal gets value \true and one gets value \false.
  
  For each clause, take a copy of the gadget and replace the
  labels \(a_1\), \(a_2\) and \(a_3\) by the literals, and
  \(\bar{a}_1\), \(\bar{a}_2\) and \(\bar{a}_3\) by the
  complements of the literals in the clause.  Add to that a
  vertex \(v_x\) for each variable \(x\), and join it to all
  vertices labeled \(x\).  Call the resulting graph \(G\), and
  let \(T\) be formed by all \(v_x\) together with the union of
  all black vertices from the gadgets.
  
  We will show a 1-1 correspondence between truth assignments for
  \(V\) that solve \(C\) and acyclic \(T\)-orientations of \(G\).
  Start with a truth assignment.  For each edge labeled
  \(x\bar{x}\), orient it from \(x\) if \(x\) is assigned \true
  and towards \(x\) otherwise.  Consider a clause and its
  respective gadget.  The three literals in the clause are not
  all \true and not all \false, so the three special edges are
  not all directed the same way.  It follows from Fact 2 that one
  can orient (uniquely) all gadgets extending these orientations.
  In this orientation, all vertices labeled with the same literal
  are sources in their gadgets if that literal is \true, and
  sinks otherwise.  This \(T\)-orientation can now be extended to
  the whole~\(G\), directing all edges incident to \(v_x\)
  towards it if \(x\) is \true and the opposite otherwise.
  
  Conversely, suppose a \(T\)-orientation of \(G\) is given.
  Since all neighbors of \(v_x\) are pairwise non-adjacent,
  \(v_x\) is either a source or a sink.  Assign \(x\)
  \true if \(v_x\) is a sink, \false otherwise.  The
  fact that each gadget is acyclically \(T\)-oriented shows that
  in the corresponding clause the not-all-equal condition is
  satisfied.
\end{proof}

Since simple powers have originated the problem in Proposition
\ref{qnpc}, we tried to look at another extension of its first
part, namely, allow only monomials with support size exactly 2.
This was short lived, though:

\begin{pro}
  The problem: 
  \begin{quote}
    \emph{given a collection \M of monomials, each with support
      of size 2, does there exist a cool ordering for \M?}
  \end{quote}
  is NP-complete.
\end{pro}
\begin{proof}
  We reduce the quadratic case to this.  Suppose that \M is a
  collection of quadratic monomials, and let \(\M'=\{xy \mid
  xy\in\M\}\cup\{x^2y \mid x^2\in\M \mbox{ and } xy\not\in\M\}\).
  It is easy to check that \M and \(\M'\) have precisely the
  same cool orderings.
\end{proof}

There is a lot of leeway in the reduction in the proof of
Proposition \ref{npc}.  For instance, the \(v_i\) could be
eliminated, and similarly labeled vertices could be merged.  One
could add irrelevant vertices of both types and show that
existence of acyclic \(T\)-orientations is NP-complete even if
\(|T|=\frac{1}{2}|VG|\) (any constant between \(0\) and \(1\)
would do).  On the extremes, the problem can be solved:

When \(T=VG\), that is recognition of comparability graphs.  When
\(T\) induces a bipartite graph, any acyclic orientation in which
one side of \(T\) consists only of sources and the other (if it
exists) only of sinks is a \(T\)-orientation.  This takes care of
\(|T|\leq 2\); actually, for any fixed \(k\), if \(|T|\leq k\),
one can restrict the search for a \(T\)-orientation to a
polynomial number of acyclic orientations that can be
systematically enumerated.

\section{Lifting the ideal}
\label{sec:pre-image}

Here we present the solution to Problem \ref{pinfg}.

\begin{teo}\label{preimage}
  Given an antichain of monomials \(M\subseteq\Xc\), the following are
  equivalent:
  \begin{itemize}
  \item [i)] \pin{\ideal{M}} is a finitely generated ideal of \Xs.
  \item [ii)] For every \(m\) in \M and any letters \(x\neq z\neq
    y\) such that \di{xy}{m}, there exists a monomial \(w\) in \M
    such that \era{w}{z} divides either \(mx^{-1}\) or \(my^{-1}\).
    (Note that \(x=y\) is included.)
  \item [iii)] For every \(m\) in \M and any letter \(z\) such that
    no power of it is in \M, if \era{m}{z} has degree \(\geq 2\),
    there exists a monomial \(z^rt\) in \(M_2\),
    such that \(t\in\su{m}\).
  \item [iv)] \pin{\ideal{M}} is generated by the inverse images of
    the monomials  \(m\in\ideal{M}\) such that, for every letter
    \(x\), \(\partial_x m\leq \max_{u\in M_2} \partial_x u\).
  \end{itemize}
\end{teo}
\begin{proof}
  (\textit{i}) implies (\textit{ii}): Given \(m\), \(x\), \(y\),
  and \(z\), choose \(u=xz^rvy\in\pinp{m}\), where \(r\geq 0\).  Now,
  for every \(s\geq r\), \(xz^svy\in\pin{\ideal{M}}\), and since
  \pin{\ideal{M}} is finitely generated, some minimal generator \(g\)
  divides infinitely many of these.  This is only possible if \(g\)
  divides some \(xz^s\) or some \(z^svy\).  So, \era{\pi(g)}{z}
  divides either \(my^{-1}\) or \(mx^{-1}\).  The result follows by
  taking any \(w\) in \M that divides~\(\pi(g)\).
  
  (\textit{ii}) implies (\textit{iii}): Given \(m\) and
  \(z\), choose letters \(x\) and \(y\) such that
  \di{xy}{\era{m}{z}}.  Let \(w\) be given by condition
  (\textit{ii}), with minimal support.  Since \M is an
  antichain, \(z\in\su{w}\).  Suppose, by contradiction, that
  \era{w}{z} has degree \(\geq 2\).  We apply condition
  (\textit{ii}) to \(w\), obtaining a \(w'\); that new monomial
  could also play the role of \(w\) with respect to \(m\), so, by
  minimality of \(w\), it cannot exist.  Since no power of \(z\)
  is in \M, \(w=z^rt\) for some letter \(t\neq z\).  As
  \era{w}{z} divides either \(mx^{-1}\) or \(my^{-1}\), it
  follows that \(t\in\su{m}\).
  
  (\textit{iii}) implies (\textit{iv}): Let \(W\) be the set
  claimed to generate \pin{\ideal{M}}.  If this is false, then
  \(\ideal{W}\subsetneq\pin{\ideal{M}}\), since
  \(W\subseteq\pin{\ideal{M}}\).  So, there must exist a
  \(w\in\pin{\ideal{M}}\), of minimum length, with no factor in
  \(W\).  So, for some letter~\(z\), \(\partial_z\pi(w)> r =
  max_{u\in M_2} \partial_z u\).
  
  Suppose that \(z^r\in M\).  Clearly \(w\) has a proper factor
  \(u\) such that \di{z^r}{\pi(u)}, so \(u\in\pin{\ideal{M}}\).
  By minimality of \(w\), \(u\) has a factor in \(W\); then, so
  does \(w\), a contradiction.  So, \(z^r\notin M\), and by the
  choice of \(r\) and as \M is an antichain, no power of \(z\)
  lies in~\M.  Now, let \(m\in M\) be such that \di{m}{\pi(w)}.
  By (\textit{iii}), there exists a \(z^st\in M\) such that
  \(t\in\su{m}\).  Since \(s\leq r\), \(w\) has a proper factor
  \(u\) such that \di{z^st}{\pi(u)}.  We get a contradiction
  again, that finishes the proof.
  
  (\textit{iv}) implies (\textit{i}): We deserve the rest.
\end{proof}

We briefly relate this result to the preceding ones.  It is easy to
check from the definitions that if \pin{\ideal{M}} is finitely
generated, then every ordering of \X is cool for \M.  This can also
be seen from the fact that if \M satisfies the condition in Theorem
\ref{preimage} (\textit{iii}), then it also satisfies the condition
of Theorem~\ref{totallycool}.  The converse is not true; the simplest
example is \(M=\{a^2,bc\}\) --- here, every ordering is cool, but
\pin{\ideal{M}} is not finitely generated.  Actually, if one starts
with any \M for which every ordering is cool and substitutes each
letter for its square, this property is preserved.  But now,
\pin{\ideal{M}} is not finitely generated.

The similarity between Theorems \ref{preimage} and \ref{FG} may
suggest that perhaps a restricted form of \ref{FG}
(\textit{ii}) involving \(M_2\) would hold.  That is not
likely, as suggested by
\(M=\{x_1x_2^2,x_1x_2x_3^2,x_1x_2x_3x_4^2,\ldots\}\); the natural
ordering of \(x_1, x_2,\ldots,x_n\) is cool for \M, for any
\(n\), even though \(M_2\) is quite skimpy.

\section{Commutative ideals given by inequalities}
\label{sec:comm-ideals-given}

Another way of giving an ideal of \Xc is as the pre-image of an ideal
under a morphism from \Xc to another free commutative monoid.  This is
useful only if there is a nice way of describing the morphisms and the
ideals of the target.  We will consider morphisms between commutative
monoids and lift ideals given by generators.

In this setting, it will be convenient to switch to an additive
notation for \Xc.  We number the letters of \X as
\(x_1,x_2,\ldots,x_n\), and identify \Xc with \(\N^n\) by the
isomorphism given by \(x_1^{i_1}x_2^{i_2}\cdots x_n^{i_n}\mapsto
\mathbf{x}=(i_1,i_2,\ldots,i_n)\).  In this notation, a set
\(I\subseteq\N^n\) is an ideal if \(\mathbf{x}\geq\mathbf{y}\in I\)
implies \(\mathbf{x}\in I\) (as usual, \(\mathbf{x}\geq\mathbf{y}\)
means \(\mathbf{x}_i\geq\mathbf{y}_i\) for every~\(i\)). Other terms
require translation: monomials become vectors, letters become indices
or coordinates of vectors, and so on.

Consider a morphism \(\varphi:\Xc\rightarrow[Y\mspace{1mu}]\),
where an isomorphism \([Y\mspace{1mu}]\rightarrow\N^m\) is
already fixed.  If \(\varphi(x_j)=\prod_i y_i^{a_{ij}}\), then
\(\varphi\) is the linear map \(\N^n\rightarrow\N^m\) given by
the matrix \(A=(a_{ij})\), where each \(a_{ij}\) is a nonnegative
integer.  If \(J\) is an ideal of \(\N^m\), then
\(I=\{\mathbf{x}\in\N^n\mid A\mathbf{x}\in J\}\) is an ideal of
\Xc; moreover, if \(J\) is generated by the finite set \(W\),
then that same ideal \(I\) can be described as:
\[
  I(A,W) = \{\mathbf{x}\in\N^n \mid A\mathbf{x}\geq w
  \,\text{for some}\, w\in W\}.
\]
When \(w\) is a vector in \(\N^m\), we write \(I(A,w)\) for
\(I(A,\{w\})\).  Also, we stress that we only consider \Iaw when
\(A\) is nonnegative.

Given a generating set \M of an ideal of \Xc, one can write down a
description \(\ideal{M} = I(I_n,M)\), simply using \(Y=X\) and the
identity morphism as \(\varphi\).  From the complexity viewpoint, we
notice that the new description has size bounded by a polynomial on
the size of \M; one interesting feature of the new type of description
is that it can be much more compact.  Just to give a trivial example,
consider, for each \(k\) in \(\N\), the ideal \(\{\mathbf{x}\in\N^2
\mid x_1+x_2\geq k\}\).  The size of this description is \(O(\log
k)\), while clearly it has \(k+1\) minimal generators.

Clearly, \(\Iaw= \bigcup_{w\in W} I(A,w)\), and this suggests the
following definition: we say that an ideal of \Xc is
\emph{convex} if it is of form \(I(A,w)\) for some integer matrix
\(A\) and vector \(w\).  The name is motivated by the following
fact, that follows from standard results in the theory of
polyhedra (see \cite{Sj} for terminology and facts about
polyhedra that we use).

\begin{pro}\label{polyhedra}
  Let \(I=\ideal{M}\) be an ideal of \(\Xc=\N^n\), with \M
  finite.  The following are equivalent:
  \begin{itemize}
  \item [i) ] \(I\) is convex.
  \item [ii) ] \(I\) is the intersection of \(\N^n\) with a
    convex set in \(\R^n\).
  \item [iii) ] \(I= \N^n \cap (\mathrm{conv}(M)+\R^n_+)\)
    (\(\mathrm{conv}(M)\) is the convex hull of \(M\)).
  \end{itemize}
\end{pro}

So, \Iaw is a union of convex ideals.  It turns out that any
union of convex ideals can be expressed as an \Iaw.  As we see
below, this can be done without wasting much space, so we can
switch descriptions without penalty in the coarse complexity of
the problems we will talk about.

\begin{lem}\label{union}
  Let \(I_1=I(A^{(1)},w^{(1)}), I_2=I(A^{(2)},w^{(2)}), \ldots,
  I_r=I(A^{(r)},w^{(r)})\) be ideals of \Xc.  Then there exist a
  matrix \(A\) and a set \(W\) of vectors, with total size
  polynomial in the total size of the descriptions
  \(I(A^{(i)},w^{(i)})\), such that \(\bigcup_i I_i = I(A,w)\).
\end{lem}
\begin{proof}
  Let \(A\) result from piling up the matrices \(A^{(i)}\) on top
  of each other.  For each \(i\), let \(w^i\) result from
  extending \(w^{(i)}\) with null entries corresponding to the
  inequalities of the other systems; so, \(I(A,w^i)=I_i\).
  Finally, let \(W=\{w^1,w^2,\ldots,w^r\}\).

  For those of a more categorical persuasion, the proof is simply the
  substitution of a family of morphisms by its direct product.
\end{proof}

Now we consider what happens to the three guiding problems of the
introduction when \(I\) is given in the form \Iaw.  Problem
\ref{existscool} is sort of hopeless: since a description
\ideal{M} can be converted into a description \Iaw of size
polynomial in the size of \M, and Problem \ref{existscool} is
NP-complete when \M is the given data.  It follows that this
problem with \(I\) given as \Iaw is NP-hard; to make things
worse, we cannot even assert that it is in NP.  At this point, we
refer the reader again to \cite{GJ} for a refresher on
NP-completeness concepts, and, in particular, to the
satisfiability problem, that will play an important role in the
remainder of this section.

For the other two problems, our results are similarly bad and
more definite. They will be shown to be coNP-complete.  Indeed,
we will add a new problem to the pack, that is completely trivial
if the ideal is given by generators:

\begin{problem}\label{Mdois}
  Given an ideal \(I\) of \Xc, is it generated by monomials with
  support of size at most \(2\)?  That is, is there a set \M such
  that \(M=M_2\) and \(\ideal{M}=I\)?
\end{problem}

We register two basic algorithms pertaining to these problems.

\begin{lem}
  Given an ideal \(I=\Iaw\) and a vector \(\mathbf{x}\), it can
  be decided in polynomial time whether \(\mathbf{x}\in I\) and
  whether \(\mathbf{x}\) is a minimal generator of \(I\).
\end{lem}
\begin{proof}
  Computing \(A\mathbf{x}\) and comparing the result with each
  member of \(W\), we quickly decide membership in \(I\).  To
  decide whether an \(\mathbf{x}\in I\) is minimal, it is enough
  to verify that each vector obtained from \(\mathbf{x}\) by
  subtracting \(1\) from a positive coordinate is not in \(I\).
\end{proof}

In what follows, the proofs will be a bit sketchy, with some bare
statements; filling in the details is routine handiwork.

\begin{pro}\label{coNP}
  Problems \ref{pinfg}, \ref{IMfg} and \ref{Mdois} are in coNP,
  when the ideal is given as \Iaw.
\end{pro}
\begin{proof}
  For each problem, when the answer to an instance is
  \textbf{no}, we will present a short certificate, verifiable in
  polynomial time.  That will be a minimal generator of the
  ideal, and some additional information.  Notice that any
  minimal generator has coordinates bounded by the maximum of all
  coordinates in members of \(W\), so it can be part of a short
  certificate.

  For Problem \ref{Mdois}, a certificate is simply a minimal
  generator with support of size at least \(3\).

  For Problem \ref{pinfg}, a certificate is a minimal generator
  \(\mathbf{m}\) and an index \(z\) such that item
  (\textit{iii}) of Theorem \ref{preimage} is violated.  That
  amounts to the following:
  \begin{itemize}
  \item [-] There is no vector in \(I\) whose support is
    \(\{z\}\) (no power of \(z\) is in \M). This happens if and
    only if, for each \(w\in W\), there is an index \(i\) such
    that \(w_i>0\) and \(a_{iz}=0\).
  \item [-] \(\sum_{i\neq z} m_i\geq 2\) (\era{m}{z} has degree
    \(\geq 2\)).
  \item [-] There is no \(\mathbf{x}\in I\) and index \(t\neq z\)
    such that \(m_t>0\), \(x_t=1\) and \(x_i=0\) for every
    \(i\neq t,z\).  This is true, for each candidate \(t\), if and
    only if for every \(w\in W\) there exists an \(i\) such that
    \(w_i>a_{it}\) and \(a_{iz}=0\).
  \end{itemize}
  
  For Problem \ref{IMfg}, we assume, without loss of generality,
  that the ordering on the letters is that of the indexing.  Now,
  a certificate consists of a minimal generator \(\mathbf{m}\)
  and an index \(x\), interior to \(\mathbf{m}\) satisfying the
  condition: there is no minimal generator \(\mathbf{s}\) whose
  first or last positive entry is in position \(x\), and such
  that if \(\mathbf{s'}\) results from \(\mathbf{s}\) by turning
  the \(x\)-component to \(0\), then \(\mathbf{s'}\leq
  \mathbf{m}\).  This condition can be checked as follows.  Let
  \(\mathbf{m_\leftarrow}\) (\(\mathbf{m_\rightarrow}\)) result
  from \(\mathbf{m}\) by changing to zero all components with
  index bigger (smaller) than \(x\).  Let also \(A'\), \(W'\)
  result from eliminating all rows \(i\) such that \(a_{ix}>0\).
  Then, \(x\) satisfies the required condition if and only if
  neither \(\mathbf{m_\leftarrow}\) nor
  \(\mathbf{m_\rightarrow}\) is in \(I(A',W')\).
\end{proof}

\begin{pro}\label{coNPc}
  Problems \ref{pinfg}, \ref{IMfg} and \ref{Mdois} are
  coNP-complete, when the ideal is given as \Iaw.
\end{pro}
\begin{proof}
  We will reduce directly from \textsc{Sat} to the negative of
  each problem.  The reductions will have a lot in common.  From
  each instance of \textsc{Sat}, we will produce a family of
  convex ideals \(I_i\), like in Lemma \ref{union}; instead of
  presenting them in matrix form, we write them as systems of
  linear inequalities.
  
  Given an instance \(S\) of \textsc{Sat} on variables
  \(x_1,x_2,\ldots,x_n\), (we assume \(n\geq 3\)) our
  inequalities will involve the variables
  \(x_1,\bar{x}_1,x_2,\bar{x}_2\ldots,x_n,\bar{x}_n\), in obvious
  correspondence to the literals.  For each clause, the
  corresponding \emph{clause inequality} is
  \[
    \text{sum of the literals in the clause}\,\, \geq 1.
  \]
    
  Let \(I_0\) be defined by the clause inequalities, together
  with the \emph{boolean inequalities} \(x_i+\bar{x}_i\geq 1\)
  for \(i=1,2,\ldots,n\).  The specific use of \(I_0\) is the
  following: \textit{\(\mathbf{x}\) is a solution of \(I_0\) in
    nonnegative integers such that each boolean inequality is
    satisfied as equality if and only if \(\mathbf{x}\) is a
    boolean assignment satisfying \(S\)}.
  
  We also define, for each \(i=1,2,\ldots,n\), the ideal \(I_i\)
  given by the single inequality \(x_i+\bar{x}_i\geq 2\).  Notice
  that \(I_i\) has three minimal generators: two with single
  support (\(x_i=2\) or \(\bar{x}_i=2\)), the other with
  two-element support (\(x_i=\bar{x}_i=1\)).
  
  \emph{Reduction to Problem \ref{Mdois}:} The instance \(P\)
  of Problem \ref{Mdois} consists precisely of the systems
  \(I_0,I_1,\ldots,I_n\).
  
  Suppose that \(S\) is satisfiable, and let \(\mathbf{x}\) be a
  boolean assignment satisfying~\(S\).  Since for each
  \(i\), exactly one of \(x_i\) or \(\bar{x}_i\) equals \(1\),
  the support of \(\mathbf{x}\) has size \(n\geq 3\), and
  \(\mathbf{x}\) is not in any \(I_i\), \(i\geq 1\).  On the
  other hand, clearly \(\mathbf{x}\) is in \(I_0\).  Also,
  \(\mathbf{x}\) is minimal, since zeroing any variable would
  violate the corresponding boolean inequality.  So, \(P\) has a
  negative answer if \(S\) is satisfiable.
  
  Conversely, suppose \(P\) has a negative answer, that is, the
  corresponding ideal has a minimal generator \(\mathbf{x}\)
  whose support has size \(\geq 3\).  Clearly it cannot be in any
  \(I_i\) with \(i\geq 1\), so it is in \(I_0\), and
  \(x_i+\bar{x}_i=1\) for each \(i\).  Hence, \(S\) is
  satisfiable.
  
  \emph{Reduction to Problem \ref{IMfg}:} We introduce two new
  variables, \(y\) and \(z\), besides the ones we already have.
  The instance \(P\) of Problem \ref{IMfg} consists of the
  systems \(I_1,\ldots,I_n\), together with \(I_0'\), which is
  \(I_0\) with the addition of the inequality \(y\geq 1\).  The
  variables of \(P\) will be ordered increasingly as
  \(y,z,x_1,\bar{x}_1,x_2,\bar{x}_2,\ldots,x_n,\bar{x}_n\).

  Suppose that \(S\) is satisfiable, define \(\mathbf{x}\) as
  before, and extend it by setting \(y=1\) and \(z=0\).  Then
  this is a minimal generator and is only in \(I_0'\).  Now,
  \(z\) is internal to this vector, but no minimal generator of
  the ideal has \(z\) in its support (let alone, as an extreme
  entry), so condition (\textit{ii}) of Theorem \ref{FG} is
  violated, and \(P\) has a negative answer.

  Conversely, if \(P\) has a negative answer, there exists a
  minimal generator and an internal variable such that condition
  (\textit{ii}) of Theorem \ref{FG} is violated.  This minimal
  generator cannot be in any of the \(I_i\), \(i\geq 1\), since
  those have no internal letters.  So it is in \(I_0'\), and must
  have \(y=1\), \(z=0\), and the other variables must be a boolean
  assignment that satisfies \(S\). The problematic internal
  variable must be \(z\), but who cares?
    
  \emph{Reduction to Problem \ref{pinfg}:} We use just one new
  variable \(y\).  The instance \(P\) of Problem \ref{pinfg}
  consists of \(I_0\), a new system \(I_*\), with the single
  inequality \(y\geq 2\), systems \(I_i'\), each obtained from
  \(I_i\) by the addition of the inequality \(y\geq 1\).  By
  arguments similar to the preceding ones and the help of Theorem
  \ref{preimage} (\textit{iii}), it can be shown that \(S\) is
  satisfiable if and only if \(P\) has a negative answer.
\end{proof}

Proposition \ref{polyhedra} (\textit{iii}) says that a convex
ideal is the set of integer points of a \emph{blocking
  polyhedron}.  Such polyhedra, and mostly their integer points,
have been the subject of a lot of attention in the context of
combinatorial and integer programming.  This, and perhaps sheer
curiosity, justify asking what happens to Problems
\ref{pinfg}--\ref{Mdois} if one restricts the questions to convex
ideals (given in the form \(I(A,w)\)).  No one of the definite
results we presented so far applies to convex ideals; in
particular, the proof of Theorem \ref{qnpc} constructs ideals
that are not convex, so even Problem \ref{existscool}'s status is
undecided.

\end{document}